\documentclass{article}

\usepackage[pdftex,
pdfauthor={U. Villa, A. Marques},
pdftitle={An UQ-ready finite element solver for a two-dimensional RANS model of free plane jets},
pdfkeywords={Finite element method, Computational Fluid Dynamics, Turbulence model},
pdfproducer={.},
pdfcreator={.}]{hyperref}

\usepackage[automark]{scrpage2}
\usepackage{authblk}
\usepackage{mathtools}
\usepackage{cite}
\ihead[]{}
\ohead[]{}
\ifoot[]{}
\ofoot[]{} 

\usepackage{color}

\newcommand{\Gammain}{\Gamma_{\rm in}}
\newcommand{\Gammaff}{\Gamma_{\rm ff}}

\title{An UQ-ready finite element solver for a two-dimensional RANS model of free plane jets\thanks{This work is supported by DARPA's EQUiPS program under contract number W911NF-15-2-0121.}}

\author[1]{Umberto Villa}

\author[2]{Alexandre Noll Marques}

\affil[1]{Institute for Computational Engineering \& Sciences, The University of Texas at Austin, USA}

\affil[2]{Department of Aeronautics \& Astronautics, Massachusetts Institute of Technology, Cambridge, MA, USA}

\usepackage{amsmath,amsthm,amsfonts,amssymb}
\usepackage{soul}
\usepackage{tabularx}
\usepackage[usenames,dvipsnames]{xcolor}
\usepackage{graphicx,hyperref,float}
\usepackage{grfext}
\usepackage{url}
\usepackage{subfigure}
\usepackage{appendix}
\usepackage{geometry}
\usepackage{algorithmic,algorithm}  
\usepackage{url}
\usepackage{tikz}
\usepackage{pgfplots}
\usepackage{pgfplotstable}
\pgfplotsset{compat=newest}

\newcommand{\bit}{\begin{itemize}}
\newcommand{\eit}{\end{itemize}}

\graphicspath{./plots}
\DeclareGraphicsExtensions{.pdf,.eps,.png,.jpg,.jpeg}

\begin{document}

\maketitle
\begin{abstract}
Numerical solution of the system of partial differential equations arising from the Reynolds-Averaged Navier-Stokes (RANS) equations with $k-\epsilon$ turbulence model presents several challenges due to the advection dominated nature of the problem and the presence of highly nonlinear reaction terms. State-of-the-art software for the numerical solution of the RANS equations address these challenges by introducing non-differentiable perturbations in the model to ensure numerical stability. However, this approach leads to difficulties in the formulation of the higher-order forward/adjoint problems, which are needed for scalable Hessian-based uncertainty quantification (UQ) methods. In this note, we present the construction of a UQ-ready flow solver, i.e., one that is smoothly differentiable and provides not only forward solver capabilities but also adjoint and higher derivatives capabilities.\\[-1mm]

\end{abstract}

\section{Introduction} \label{sec:intro}

Free turbulent jets are prototypical flows believed to represent the dynamics in many engineering applications, such as combustion and propulsion. As such, free jet flows are the subject of several experimental \cite{Gutmark76_planeJet, Gutmark78_impJet, Krothapalli81_mixing}
and numerical investigations \cite{Zhou99_RANS, Ribault99, Stanley02_mixing, Klein03_DNSplanejet, klein15} and constitute an important benchmark for turbulent flows.
Here we restrict our attention to non-reactive plane jets modeled by the two-dimensional Reynolds-Averaged Navier-Stokes (RANS) equations and the $k-\epsilon$ turbulence model.
This document describes in details the formulation and discretization used to create a finite element solver of such flows.
The main goal is to construct
a \emph{UQ-ready} flow solver, i.e., one that is smoothly
differentiable and provides not only forward solver capabilities but
also adjoint and higher derivatives capabilities.
This solver is implemented in FEniCS~\cite{AlnaesBlechta2015a}.
In Section~\ref{sec:model} we present the model used to represent the physics of the jet.
Next, in Section~\ref{sec:discretization} we discuss the discretization of the physical model and in Section~\ref{sec:validation} we present 
numerical validation results.

\section{Free plane jet model} \label{sec:model}
We consider a free plane jet in conditions similar to the ones reported in \cite{Klein03_DNSplanejet, klein15}.
Namely, the flow exits a rectangular nozzle into quiescent surroundings with a prescribed top-hat velocity profile and turbulence intensity.
The nozzle has width $D$, and is infinite along the span-wise direction.

\begin{figure}[h!]
	\begin{center}
		\includegraphics[width=3.5in]{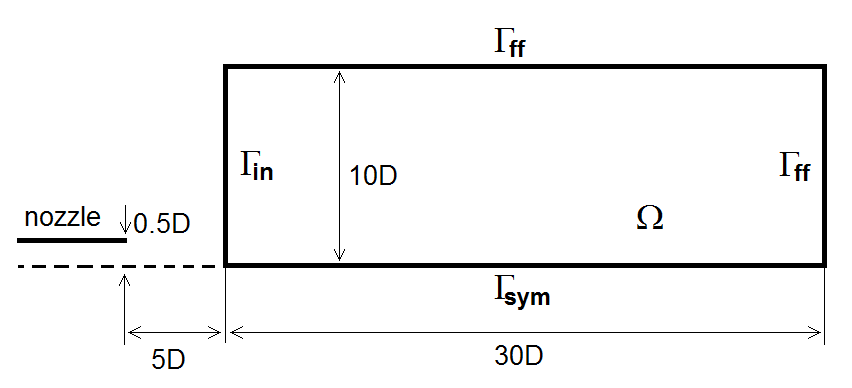}
	\end{center}
	\caption{Illustration of the free plane jet setup.
		The diameter of the nozzle is denoted by $D$. The simulation domain $\Omega$ is composed of a $30D \times 10D$ box situated at a distance $5D$ downstream to the nozzle exit.}
	\label{fig:setup}
\end{figure}

Our simulation model computes the flow in a rectangular domain $\Omega$ located at a distance $5D$ downstream from the exit of the jet nozzle, as illustrated in Figure \ref{fig:setup}.
By doing so, modeling the conditions at the exit plane of the jet nozzle is avoided.
Instead, direct numerical simulation data are used to define inlet conditions at the surface $\Gammain$.
The dynamics are modeled with the steady incompressible Reynolds-Averaged Navier-Stokes equations, complemented by the $k-\epsilon$ turbulence model \cite{Launder74_ke}:
\begin{align}
R_v &:= (\boldsymbol{v} \cdot \nabla)\boldsymbol{v} 
+ \dfrac{1}{\rho}\nabla p 
- \nabla \cdot ((\nu + \nu_t)\bar{\bar{S}}(\boldsymbol{v}))
= 0,
\label{eq:mom}\\
R_p &:= \nabla \cdot \boldsymbol{v} = 0,
\label{eq:mass}\\
R_k &:= \boldsymbol{v} \cdot \nabla k 
- 2\nu_t(\bar{\bar{S}}(\boldsymbol{v}):\bar{\bar{S}}(\boldsymbol{v}))
+ \epsilon 
- \nabla \cdot \left(\left(\nu + \dfrac{\nu_t}{\sigma_k}\right)\nabla k \right) = 0,
\label{eq:k}\\
R_{\epsilon} &:= \boldsymbol{v} \cdot \nabla \epsilon
- 2C_{1\epsilon}\dfrac{\epsilon\nu_t}{k}(\bar{\bar{S}}(\boldsymbol{v}):\bar{\bar{S}}(\boldsymbol{v})) 
+ C_{2\epsilon}\dfrac{\epsilon^2}{k} 
- \nabla \cdot \left(\left(\nu + \dfrac{\nu_t}{\sigma_{\epsilon}}\right) \nabla \epsilon \right) = 0,
\label{eq:e}
\end{align}
where $\boldsymbol{v} = [v_x, v_y]$ denotes the velocity vector, $p$ denotes pressure, $\rho$ is the density, $\nu$ is the kinematic viscosity, and $\bar{\bar{S}}$ is the strain rate tensor given by
\begin{equation*}
\bar{\bar{S}}(\boldsymbol{v}) = \dfrac{1}{2}(\nabla \boldsymbol{v} + (\nabla \boldsymbol{v})^T).
\end{equation*}
In the $k-\epsilon$ turbulence model, $k$ denotes the turbulent kinetic energy, $\epsilon$ denotes the turbulent dissipation, $\nu_t$ denotes the turbulent kinematic viscosity,
\begin{equation*}
\nu_t = C_{\mu}\dfrac{k^2}{\epsilon},
\end{equation*}
and the constants of the model are
\begin{align*}
C_{\mu} = 0.09 && 
\sigma_k = 1.00, && 
\sigma_{\epsilon} = 1.30, && 
C_{1\epsilon} = 1.44, && \text{and} &&
C_{2\epsilon} = 1.92.
\end{align*}
Furthermore, the double dot product for a rank two tensor is defined as $\bar{\bar{A}}:\bar{\bar{B}} = \sum_i\sum_jA_{ij}B_{ij}$.
In the remainder of this document we refer to (\ref{eq:mom}--\ref{eq:e}) in an abbreviated form by introducing the residual operator $\boldsymbol{R} = [R_v, R_p, R_k, R_{\epsilon}]$, which acts on the state variable $\boldsymbol{s} = [\boldsymbol{v}, p, k, \epsilon]$.
Then, the equations that govern the flow dynamics are represented by
\begin{equation}\label{eq:jet}
\boldsymbol{R}(\boldsymbol{s}) = \boldsymbol{0}.
\end{equation}

These governing equations are augmented by appropriate inflow and outflow boundary conditions.
At the inlet surface $\Gammain$ Dirichlet boundary conditions are imposed on the velocity field and turbulent variables. Specifically, we set
\begin{align}\label{eq:inflow_bc}
 \boldsymbol{v}|_{\Gammain} &= \boldsymbol{v_{\text{ref}}}, &
 k|_{\Gammain} &= k_{\text{ref}}, &
 \epsilon|_{\Gammain} = \epsilon_{\text{ref}},
\end{align}
where $v_{\text{ref}}$, $k_{\text{ref}}$, and $\epsilon_{\text{ref}}$ are reference profiles obtained from DNS data or analytical approximate solutions.

At the symmetry axis surface, $\Gamma_{\text{sym}}$, no-flux boundary conditions are imposed through a combination of Dirichlet and Neumann conditions of the form
\begin{align*}
 v_y|_{\Gamma_{\text{sym}}} &= 0, &
 \left.
 \dfrac{\partial v_x}{\partial n}
 \right|_{\Gamma_{\text{sym}}} &= 0, &
 \left.
 \dfrac{\partial k}{\partial n}
 \right|_{\Gamma_{\text{sym}}} &= 0, &
 \left.
 \dfrac{\partial \epsilon}{\partial n}
 \right|_{\Gamma_{\text{sym}}} &= 0.
\end{align*}

Finally, at the  surface $\Gammaff$ we impose \emph{far-field} conditions that allow the entrainment of air around the jet.
Let $v_t$ and $v_n$ denote the tangential and normal components of
the velocity vector on $\Gammaff$.
Furthermore, let $\Gammaff^-$ be a subset of $\Gammaff$ on which $\boldsymbol{v} \cdot \boldsymbol{n} < 0$ (inflow), and $\Gammaff^+ = \Gammaff - \Gammaff^-$ (outflow).
Then, the boundary conditions at $\Gammaff$ are given by
\begin{align}\label{eq:ff}
v_t|_{\Gammaff} &= 0, &
\left.
\dfrac{\partial v_n}{\partial n}
\right|_{\Gammaff} &= 0, &
k|_{\Gammaff^-} &= 0, &
\left.
\dfrac{\partial k}{\partial n}
\right|_{\Gammaff^+} &= 0, &
\epsilon|_{\Gammaff^-} &= 0, &
\left.
\dfrac{\partial \epsilon}{\partial n}
\right|_{\Gammaff^+} &= 0.
\end{align}
These conditions imply that the flow is orthogonal to $\Gammaff$ and enters the domain as laminar, and that changes normal to the boundary are negligible. In addition, to guarantee well posedness of the advection dominated equations for the turbulence variables $k$ and $\epsilon$, we impose a homogeneus Dirichlet condition on $k$ and $\epsilon$ everywhere the flow enters the computational domain $\Omega$ through $\Gammaff^-$ and homogeneus Newmann condition everywhere the flow exits the computational domain $\Omega$ through $\Gammaff^+$.
Although we expect these conditions to hold for $\Gammaff$ sufficiently far from the jet nozzle, these conditions are only approximately valid for a finite computational domain, such as the one considered here.
We mitigate this issue by imposing \eqref{eq:ff} weakly, as suggested by \cite{Bazilevs07_weakDirichlet} and described in Section~\ref{sub:ff}.

\section{Numerical discretization}
\label{sec:discretization}

The model equations described above are solved numerically using a finite element discretization implemented on FEniCS~\cite{AlnaesBlechta2015a}.
In Section~\ref{sub:spaces} we introduce the function spaces used in the weak formulation of the free plane jet model and in its discretization.
Then, in Section~\ref{sub:weak} we present the weak formulation itself.
Next, in Section~\ref{sub:ff} we discuss the treatment of boundary conditions at the far-field.
Finally, in Section~\ref{sub:discrete} we present the final form of the discrete equations.

\subsection{Function spaces and norms}
\label{sub:spaces}

In this section we introduce the function spaces used in the remainder of this document.
First, let us denote the space of square integrable scalar functions over the domain $\mathcal{D} \subset \mathbb{R}^d$ by $\mathcal{L}_2(\mathcal{D})$:
\begin{equation*}
 \mathcal{L}_2(\mathcal{D}) := \left\{ u \, \middle| \, \int_{\mathcal{D}} u^2 \, d\mathbf{x} < \infty \right\}.
\end{equation*}
This space is equipped with the standard inner product and norm
\begin{align}
 \langle u, v \rangle_\mathcal{D} &= \int_{\mathcal{D}} u\, v \, d\mathbf{x}, \label{eq:inner}\\
 || u ||_{\mathcal{L}_2(\mathcal{D})} &= \langle u, u \rangle_\mathcal{D}^{1/2}.
\end{align}
Similarly, we will denote with  $\boldsymbol{\mathcal{L}}_2(\mathcal{D})$ and $\bar{\bar{\boldsymbol{\mathcal{L}}}}_2(\mathcal{D})$ the vectorial and tensorial counterpart of $\mathcal{L}_2(\mathcal{D})$, respectively. Then for vector functions $\boldsymbol{u}, \boldsymbol{v} \in \boldsymbol{\mathcal{L}}_2(\mathcal{D})$ and tensor functions $\bar{\bar{S}}, \bar{\bar{T}} \in \bar{\bar{\boldsymbol{\mathcal{L}}}}_2(\mathcal{D})$, we define the inner products
\begin{equation}\label{eq:inner2}
 \langle \boldsymbol{u}, \boldsymbol{v} \rangle_\mathcal{D} = \int_{\mathcal{D}}  \boldsymbol{u}\cdot \boldsymbol{v} \, d\mathbf{x}, \quad
 \langle \bar{\bar{S}}, \bar{\bar{T}} \rangle_\mathcal{D} = \int_{\mathcal{D}}  \bar{\bar{S}}:\bar{\bar{T}} \, d\mathbf{x}
\end{equation}
In what follow, we will write $\langle \cdot, \cdot \rangle$ when $\mathcal{D} = \Omega$.

Finally, we denote the Sobolev space of square integrable scalar  functions and square integral derivatives in $\Omega$ by $\mathcal{H}^1(\Omega)$:
\begin{equation*}
 \mathcal{H}^1(\Omega) := \left\{ u \in \mathcal{L}_2(\Omega) \, \middle| \, \nabla u \in \boldsymbol{\mathcal{L}}_2(\Omega) \right\},
\end{equation*}
and similarly, for vectorial functions,
\begin{equation*}
 \boldsymbol{\mathcal{H}}^1(\Omega) := \left\{ \boldsymbol{u} \in \boldsymbol{\mathcal{L}}_2(\Omega) \, \middle| \, \nabla \boldsymbol{u} \in \bar{\bar{\mathcal{L}}}_2(\Omega) \right\}.
\end{equation*}


\subsection{Weak formulation}
\label{sub:weak}

A finite element discretization of \eqref{eq:jet} requires the equations to be cast in a weak form.
We first introduce the appropriate function spaces for the solution and test functions:
\begin{align*}
& \text{\emph{solution spaces}}
\\
\mathcal{V} &:= \{\boldsymbol{v} \in \boldsymbol{\mathcal{H}}^1(\Omega) \mid \boldsymbol{v} = \boldsymbol{v_{\text{ref}}} \text{ on } \Gammain \text{ and } v_y = 0 \text{ on } \Gamma_{\text{sym}}\},
\\
\mathcal{Q} &:= \mathcal{L}_2(\Omega),
\\
\mathcal{K} &:= \{k \in \mathcal{H}^1(\Omega) \mid k = k_{\text{ref}} \text{ on } \Gammain\},
\\
\mathcal{E} &:= \{\epsilon \in \mathcal{H}^1(\Omega) \mid \epsilon = \epsilon_{\text{ref}} \text{ on } \Gammain\},
\\
\mathcal{S} &:= \mathcal{V} \times \mathcal{Q} \times \mathcal{K} \times \mathcal{E},
\\[5pt]
& \text{\emph{test spaces}}
\\
\mathcal{W} &:= \{\boldsymbol{w} \in \boldsymbol{\mathcal{H}}^1(\Omega) \mid \boldsymbol{w} = \boldsymbol{0}\ \text{ on } \Gammain  \text{ and } w_y = 0 \text{ on } \Gamma_{\text{sym}}\},
\\
\mathcal{Z} &:= \{z \in \mathcal{H}^1(\Omega) \mid z = 0 \text{ on } \Gammain\},
\\
\mathcal{Y} &:= \mathcal{W} \times \mathcal{Q} \times \mathcal{Z} \times \mathcal{Z}. 
\end{align*}

Then, the weak formulation of \eqref{eq:jet} becomes
\begin{equation}
\label{eq:weak}
\text{find } \boldsymbol{s} \in \mathcal{S} \text{ such that }
\boldsymbol{\mathcal{R}}(\boldsymbol{s}, \boldsymbol{s_t}) = \boldsymbol{0}
\text{ for all } \boldsymbol{s_t} \in \mathcal{Y},
\end{equation}
where $\boldsymbol{s_t} = [\boldsymbol{w}, q, r, u]$, $\boldsymbol{\mathcal{R}} = [\mathcal{R}_v, \mathcal{R}_p, \mathcal{R}_k, \mathcal{R}_{\epsilon}]$, and
\begin{align*}
\mathcal{R}_v &:= \left\langle (\boldsymbol{v} \cdot \nabla)\boldsymbol{v}, \boldsymbol{w} \right\rangle
-\left\langle p, \nabla \cdot \boldsymbol{w} \right\rangle
+2\left\langle (\nu + \widetilde{\nu_t}) \bar{\bar{S}}(\boldsymbol{v}), \bar{\bar{S}}(\boldsymbol{w}) \right\rangle
+ \mathcal{R}_v^{\text{ff}} = 0,
\\
\mathcal{R}_p &:= \left\langle \nabla \cdot \boldsymbol{v}, q \right\rangle = 0,
\\
\mathcal{R}_k &:= \left\langle \boldsymbol{v} \cdot \nabla k + \epsilon - 2\widetilde{\nu_t}\,(\bar{\bar{S}}(\boldsymbol{v}):\bar{\bar{S}}(\boldsymbol{v})), r \right\rangle 
+ \left\langle \left(\nu + \dfrac{\widetilde{\nu_t}}{\sigma_k}\right) \nabla k, \nabla r \right\rangle
+ R_k^{\text{ff}} = 0,
\\ 
\mathcal{R}_{\epsilon} &:= \left\langle \boldsymbol{v} \cdot \nabla\epsilon 
- 2C_{1\epsilon}\widetilde{\gamma}\widetilde{\nu_t}(\bar{\bar{S}}(\boldsymbol{v}):\bar{\bar{S}}(\boldsymbol{v})) + C_{2\epsilon}\widetilde{\gamma}\epsilon, u \right\rangle 
+ \left\langle \left(\nu + \dfrac{\widetilde{\nu_t}} {\sigma_{\epsilon}}\right) \nabla \epsilon, \nabla u \right\rangle + \mathcal{R}_{\epsilon}^{\text{ff}} = 0,
\end{align*}
where $\widetilde{\nu_t} = C_\mu \frac{k^2}{\epsilon^+}$ and $\widetilde{\gamma} = \frac{\epsilon}{k^+}$ with $\epsilon^+$ and $k^+$ being appropriately-mollified (and thus smoothly differentiable) positive functions described in Section \ref{sub:positivity}.
In the equations above, $\langle .\,,\,. \rangle$ denotes the inner product defined in \eqref{eq:inner} and \eqref{eq:inner2}, and $\mathcal{R}^{\text{ff}}_{(.)}$ correspond to the contributions due to the far-field boundary conditions, which are described below.

\subsection{Far-field boundary conditions}
\label{sub:ff}

On the far-field boundary we impose the Dirichlet conditions \eqref{eq:ff} in a weak sense. Weak imposition of Dirichlet boundary conditions, in fact, allows for more accurate boundary layer solutions of both advection-diffusion and incompressible Navier-Stokes equations, as shown in \cite{Bazilevs07_weakDirichlet}. In fact, the boundary conditions for the velocity field $\boldsymbol{v}$ are only expected to hold for $\Gammaff$ far enough from the jet nozzle and may lead to unphysical oscillations if enforced strongly in a finite computational domain. In addition, for the turbulent variables $k$ and $\epsilon$ we need a differentiable way to switch between Dirichlet and Neumann boundary condition: The region $\Gammaff^{-} \subset \Gammaff$, where the stability of those advection-dominated equations for $k$ and $\epsilon$ requires imposing Dirichlet boundary conditions, is not known \emph{a priori} but depends on the velocity field $\boldsymbol{v}$. 

Since the weak formulation naturally incorporates the Neumann conditions, to properly enforce the weak Dirichlet conditions, we will then augment  \eqref{eq:weak} with the terms
\begin{align*}
\mathcal{R}_v^{\text{ff}} &=
\left\langle \dfrac{C_b^I}{h} (\nu + \widetilde{\nu_t}) v_t, w_t \right\rangle_{\Gammaff}
- \left\langle (\nu + \widetilde{\nu_t}) (\bar{\bar{S}}(\boldsymbol{v}) \cdot \boldsymbol{n})\cdot \boldsymbol{t}, w_t \right\rangle_{\Gammaff} 
- \left\langle (\nu + \widetilde{\nu_t}) (\bar{\bar{S}}(\boldsymbol{w}) \cdot \boldsymbol{n})\cdot \boldsymbol{t}, v_t \right\rangle_{\Gammaff}, \\ 
\mathcal{R}_k^{\text{ff}} &=
\left\langle \chi(v_n) \dfrac{C_b^I}{h} \left(\nu + \dfrac{\widetilde{\nu_t}}{\sigma_k}\right) k, r \right\rangle_{\Gammaff} \\
&- \left\langle \chi(v_n) \left(\nu + \dfrac{\widetilde{\nu_t}}{\sigma_k}\right) \nabla k \cdot \boldsymbol{n}, r \right\rangle_{\Gammaff} 
- \left\langle \chi(v_n) \left(\nu + \dfrac{\widetilde{\nu_t}}{\sigma_k}\right) \nabla r \cdot \boldsymbol{n}, k \right\rangle_{\Gammaff}, \\ 
\mathcal{R}_{\epsilon}^{\text{ff}} &=
\left\langle \chi(v_n) \dfrac{C_b^I}{h} \left(\nu + \dfrac{\widetilde{\nu_t}}{\sigma_{\epsilon}}\right) \epsilon, u \right\rangle_{\Gammaff}\\
&- \left\langle \chi(v_n) \left(\nu + \dfrac{\widetilde{\nu_t}}{\sigma_k}\right) \nabla \epsilon \cdot \boldsymbol{n}, u \right\rangle_{\Gammaff}
- \left\langle \chi(v_n) \left(\nu + \dfrac{\widetilde{\nu_t}}{\sigma_k}\right) \nabla u \cdot \boldsymbol{n}, \epsilon \right\rangle_{\Gammaff},
\end{align*}
where $\chi$ denotes a smooth approximation to the indicator function for $\Gammaff^-$:
\begin{equation*}
\chi(v_n) = \dfrac{1}{2}\left(1 - \tanh\left(\dfrac{v_n}{\nu}\right)\right).
\end{equation*}

To guarantee stability, we need $C_b^I$ and $h$ such that
\begin{equation*}
||\nabla \boldsymbol{w}||_{\mathcal{L}_2(\Gammaff)} \le
\dfrac{C_b^I}{2h}||\boldsymbol{w}||_{\mathcal{L}_2(\Gammaff)}
\quad \forall \boldsymbol{w} \in \mathcal{W}.
\end{equation*}
In the discrete equations we replace $h$ by $h^\tau$, which denotes a measure of the size of each boundary element $\tau$. Finally, we choose $C_b^I = 1 \times 10^5$.

\subsection{Positivity of $k$-$\epsilon$}
\label{sub:positivity}

The most delicate issue in the solution of the RANS model is the
possible loss of positivity of the turbulence
variables. State-of-the-art software for the numerical solution of the
RANS equations (such as Featflow and OpenFOAM) address this issue by
introducing non-differentiable perturbations in the model to prevent
the solutions from becoming negative,
see e.g. \cite{lew2001note, kuzmin2007implementation}.
However, this approach is not
desirable in our case since it leads to difficulties in the
formulation of the higher-order forward/adjoint problems, which are
needed in our scalable Hessian-based UQ methods. To avoid
these issues, we introduced an appropriately-mollified (and thus
smoothly differentiable) \emph{max} function to ensure positivity of
$k$ and $\epsilon$. Namely, for a given $\varepsilon \gg 1$, we define
the positive functions
\begin{equation}
k^+ = \frac{k + \sqrt{k^2+\varepsilon^2}}{2}, \quad \epsilon^+ = \frac{\epsilon + \sqrt{\epsilon^2+\varepsilon^2}}{2}.
\end{equation}
Note that $k^+, \epsilon^+$ only appears in the definition of $\widetilde{\nu_t}$ and $\widetilde{\gamma}$, thus the nodal values of
$k$ and $\epsilon$ are not modified directly.

An alternative approach is to solve for the logarithms of $k$ and $\epsilon$ \cite{ilinca1998positivity}. However, we preferred not to follow this
approach due to the additional convective term  that appears in the modified set of equations involving exponential of the unknowns.

\subsection{Discrete equations}
\label{sub:discrete}

The discrete equations are obtained by representing the solution and test functions in finite-dimensional function spaces, and enforcing the weak statement \eqref{eq:weak} over these spaces.
Consider a triangulation of the domain into $n_{\text{el}}$ elements, denoted by $\mathcal{T}^h(\Omega)$, and let $\Omega^{e}$ denote the subdomain associated with element $e$, $e = 1,\ldots,n_{\text{el}}$ .
Then, we represent the solution and the test functions using the following finite-dimensional spaces:
\begin{align*}
& \text{\emph{solution spaces}}
\\
\mathcal{V}^h &:= \{\boldsymbol{v} \in \mathcal{V} \mid \boldsymbol{v}|_{\Omega^e} \in \mathcal{P}_2(\Omega^e) \times \mathcal{P}_2(\Omega^e) \},
\\
\mathcal{Q}^h &:= \{q \in \mathcal{Q} \mid q|_{\Omega^e} \in \mathcal{P}_1(\Omega^e)\},
\\
\mathcal{K}^h &:= \{k \in \mathcal{K} \mid k|_{\Omega^e} \in \mathcal{P}_1(\Omega^e) \},
\\
\mathcal{E}^h &:= \{\epsilon \in \mathcal{E} \mid \epsilon|_{\Omega^e} \in \mathcal{P}_1(\Omega^e) \},
\\
\mathcal{S}^h &:= \mathcal{V}^h \times \mathcal{Q}^h \times \mathcal{K}^h \times \mathcal{E}^h,
\\[5pt]
& \text{\emph{test spaces}}
\\
\mathcal{W}^h &:= \{\boldsymbol{w} \in \mathcal{W} \mid \boldsymbol{w}|_{\Omega^e} \in \mathcal{P}_2(\Omega^e) \times \mathcal{P}_2(\Omega^e) \},
\\
\mathcal{Z}^h &:= \{z \in \mathcal{Z} \mid z|_{\Omega^e} \in \mathcal{P}_1(\Omega^e) \},
\\
\mathcal{Y}^h &:= \mathcal{W}^h \times \mathcal{Q}^h \times \mathcal{Z}^h \times \mathcal{Z}^h,
\end{align*}
where $\mathcal{P}_m(\Omega^e)$ denotes the space of a polynomial functions of degree $m$ on $\Omega^e$.
This choice of function spaces corresponds to using the standard Taylor-Hood elements \cite{Donea03_FEM} for $\boldsymbol{v}$ and $p$, augmented by linear representations of $k$ and $\epsilon$.

In addition, we employ a strongly consistent self-adjoint numerical
stabilization technique (Galerkin Least Squares stabilization, \cite{Franca89_GLS, Hughes07_stab}) to
address the convection dominated nature of the RANS equations.
Specifically, we add to the discrete equations the stabilization term

\begin{equation}
\label{eq:stab}
\boldsymbol{\mathcal{R}^{\text{stab}}}(\mathbf{s^h}, \mathbf{s_t^h}) =
\sum_{e=1}^{n_{\text{el}}}
\left\langle \left. \dfrac{\partial \boldsymbol{R}}{\partial \mathbf{s^h}}\right|_{\mathbf{s^h}}\mathbf{s_t^h},
\, \bar{\bar{\tau}}^\text{e}\, \boldsymbol{R}(\mathbf{s^h}) \right\rangle_{\Omega_e},
\end{equation}
where
\begin{equation*}
\bar{\bar{\tau}}^e =
\begin{bmatrix}
\tau_v^e\\
& \tau_p^e\\
&& \tau_k^e\\
&&& \tau_{\epsilon}^e
\end{bmatrix},
\end{equation*}
and
\begin{align*}
\tau_p^e &= (h^e)^2||\boldsymbol{v}||_{\Omega^e},\\
\tau_i^e &= \dfrac{h^e}{2||\boldsymbol{v}||_{\Omega^e}}\left(\coth(\text{Pe}_i) - \dfrac{1}{\text{Pe}_i}\right), \qquad i = v,\, k,\, \epsilon.
\end{align*}
In the equations above, $\text{Pe}_i$ denotes the effective \emph{P\'{e}clet} number for each of the equations that involve convective and diffusive terms:
\begin{align*}
\text{Pe}_v &=
\dfrac{||\boldsymbol{v}||_{\Omega^e}h^e}{2(\nu + \nu_t)}, &
\text{Pe}_k &=
\dfrac{||\boldsymbol{v}||_{\Omega^e}h^e}{2(\nu + \nu_t/\sigma_k)}, &
\text{Pe}_{\epsilon} &=
\dfrac{||\boldsymbol{v}||_{\Omega^e}h^e}{2(\nu + \nu_t/\sigma_{\epsilon})},
\end{align*}
and $h^e$ denotes a measure of the size of element $e$.
Furthermore, the Jacobian of the strong residual in \eqref{eq:stab} is computed using FEniCS's symbolic differentiation capabilities.

Solving the discrete equations then amounts to finding $\mathbf{s^h} \in \mathcal{S}^h$ such that 
\begin{equation} \label{eq:discrete}
\boldsymbol{\mathcal{R}}(\mathbf{s^h}, \mathbf{s_t^h}) + \boldsymbol{\mathcal{R}^{\text{stab}}}(\mathbf{s^h}, \mathbf{s_t^h}) = \mathbf{0},
\end{equation}
for all $\mathbf{s_t^h} \in \mathcal{Y}^h$.

To solve the nonlinear system of equations \eqref{eq:discrete} that arise from
the finite element discretization of the steady-state RANS equations,
we employed a damped Newton method. The finite element
discretization is implemented in FEniCS: the weak form of the residual
(which includes GLS stabilization and mollified versions of the
positivity constraints on $k$ and $\varepsilon$ and the switching
boundary condition on the outflow boundary) was specified and we used
symbolic differentiation to obtain the expressions for the bilinear
form of the state Jacobian operator. For
robustness and global converge of the Newton method, we used
psuedo-time continuation, guaranteeing global convergence to a
physically stable solution.
Specifically, an initial guess for the solution of these equations is computed by marching the system in pseudo-time using a backward Euler discretization. Each time step requires the solution to a non-linear system, and a standard Newton iteration is used.
Then, when the pseudo-time stepping achieves $||\boldsymbol{R}|| < 0.001$, the steady equations \eqref{eq:discrete} are solved, again using Newton iteration.

\section{Validation}\label{sec:validation}
Let $D$ denote the diameter of the nozzle.
The computational domain $\Omega$ is a rectangle with lower left corner at $(5D,0)$  and upper right corner at $(35D,10D)$.
At the inlet surface $\Gammain$ Dirichlet boundary conditions \eqref{eq:inflow_bc} are imposed from the direct numerical simulation data, thus avoiding the need to model flow condition at the exit of the nozzle.
Specifically, the direct numerical simulation data described in \cite{klein15} are used to determine reference inlet profiles for velocity, $\boldsymbol{v_{\text{ref}}}$, and for turbulent kinetic energy, $k_{\text{ref}}$. The turbulent dissipation $\epsilon_{\text{ref}}$ at the inlet is then estimated by assuming a mixing length model,
$$\epsilon_{\text{ref}} = C_{\mu}\dfrac{k^{3/2}}{\ell_m}, $$ 
where $\ell_m = 0.1D$ denotes the \emph{assumed} mixing length.

\begin{figure}
\includegraphics[width=0.45\textwidth]{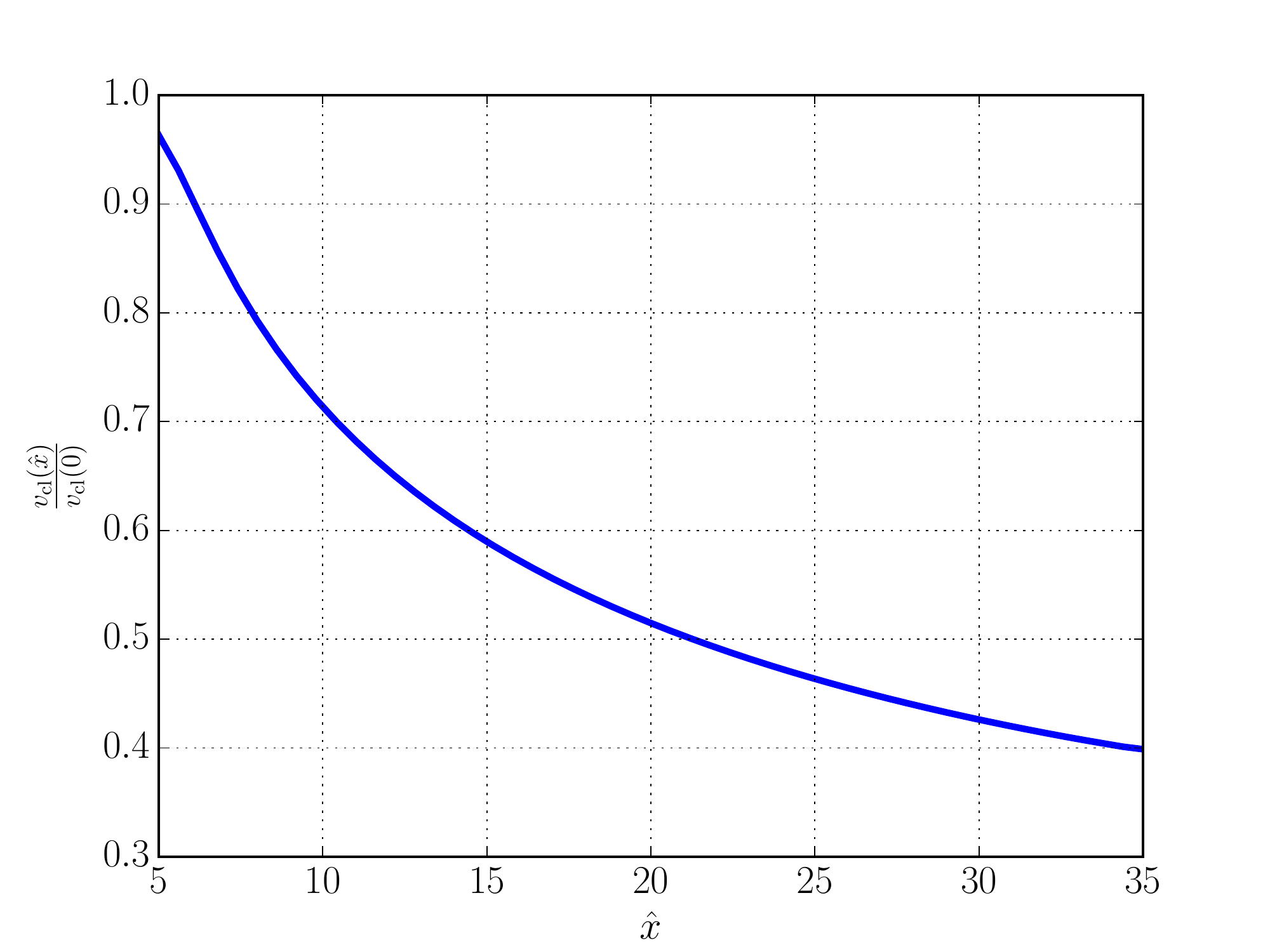}
\includegraphics[width=0.45\textwidth]{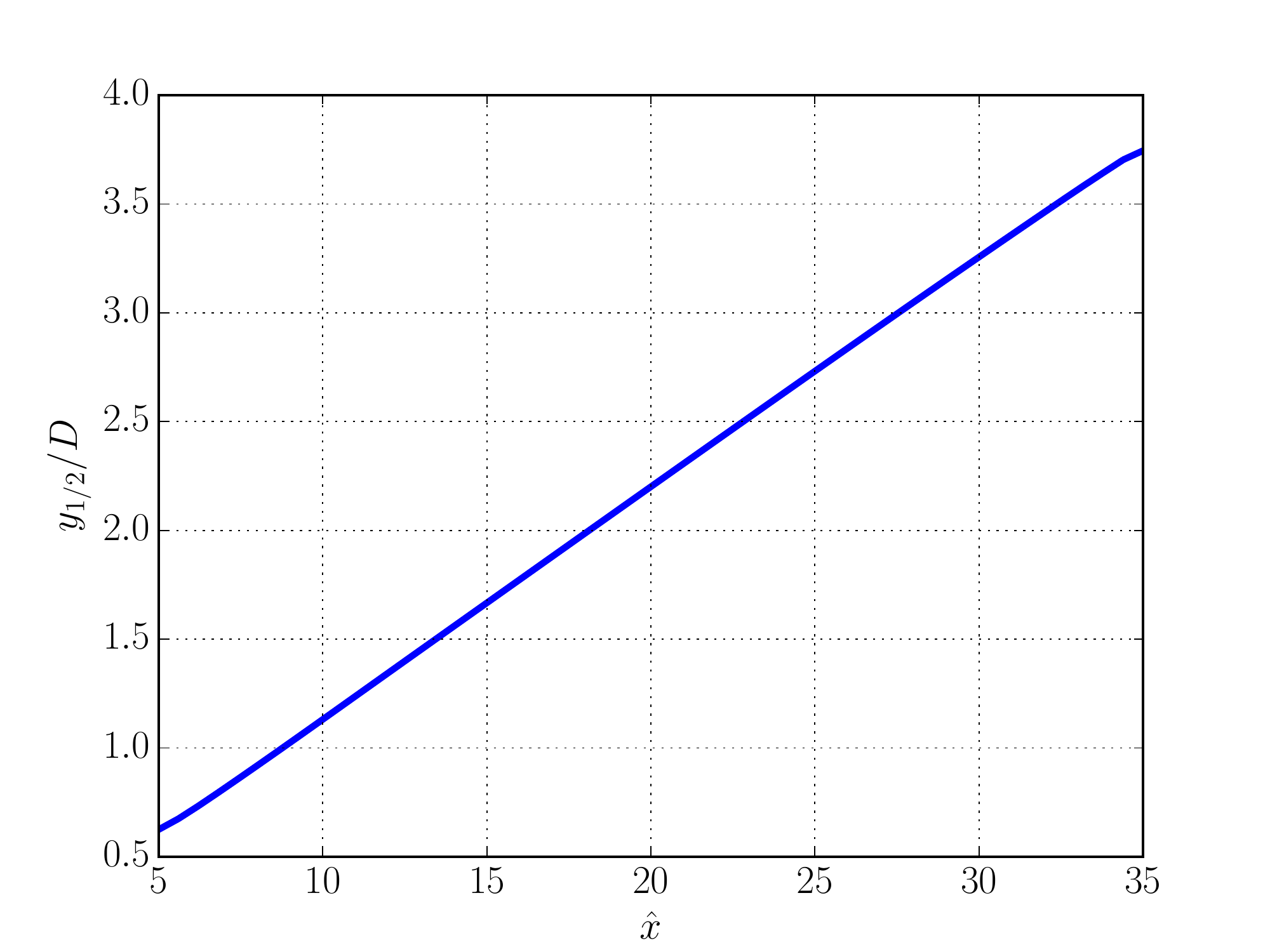}\\
\includegraphics[width=0.45\textwidth]{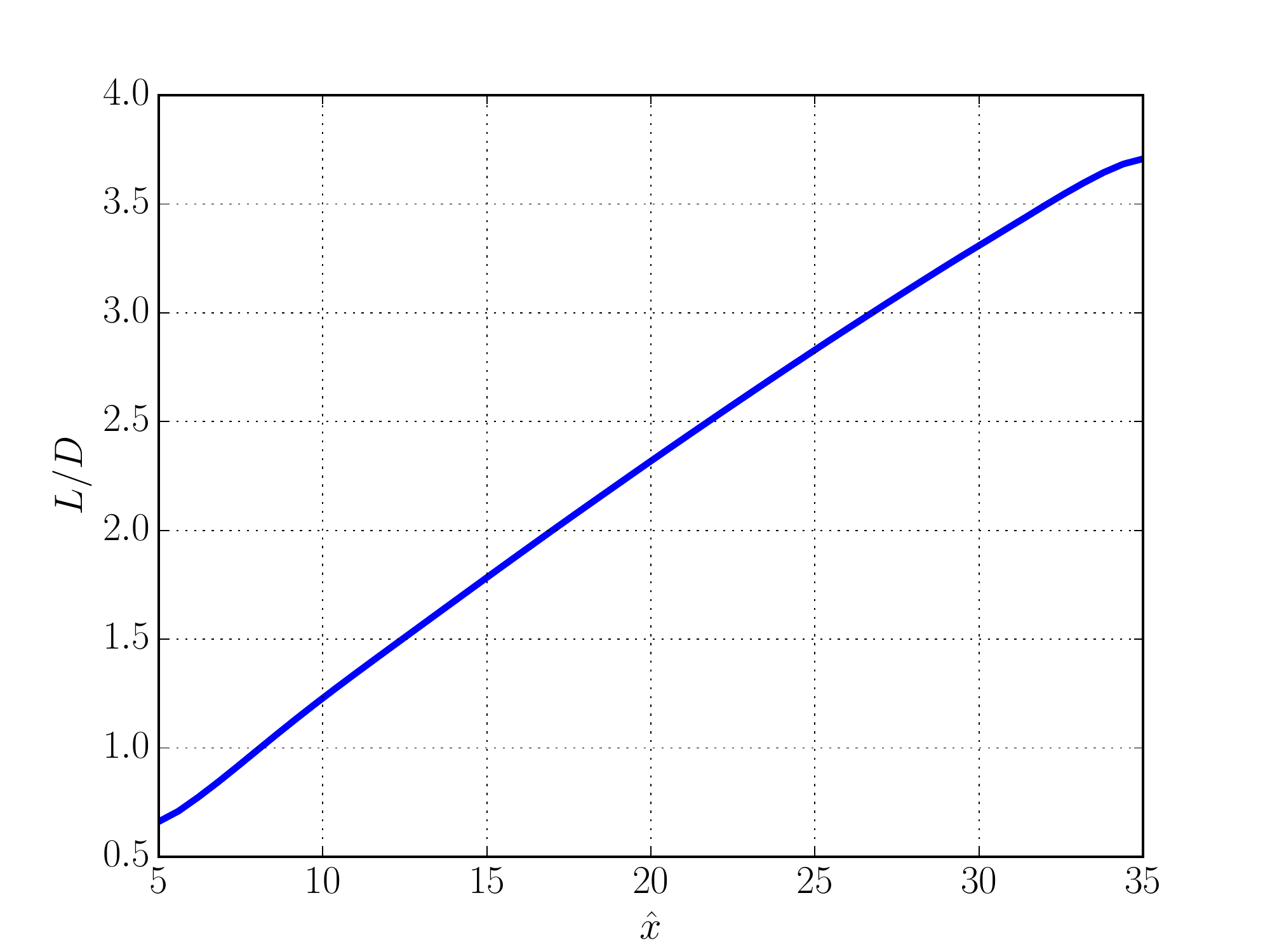}
\includegraphics[width=0.45\textwidth]{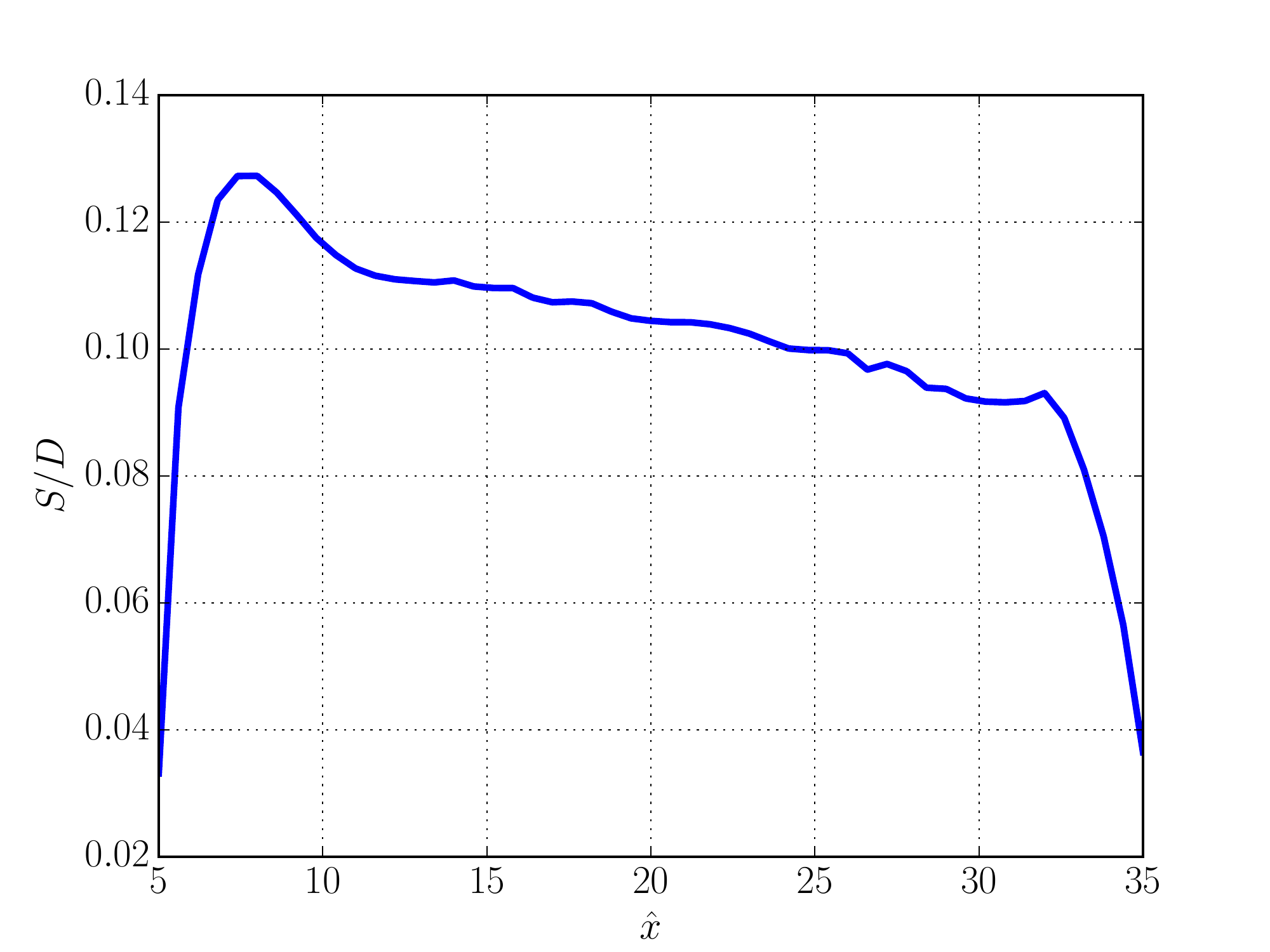}
\caption{Centerline velocity (top left), jet width $y_{1/2}$ (top right), integral jet width $L$ (bottom left), spread $S$ (bottom right).}
\label{fig:u_cl}
\end{figure}

In Figure \ref{fig:u_cl}, we show the centerline velocity $v_{\rm cl}(x) = v_x(x, y=0)$, the jet width $y_{1/2}$ (i.e. the $y$ coordinate such that $v_x(x, y_{1/2}) = v_{\rm cl}(x)$ for $x \in [5D, 35D]$), the integral jet width
$$ L(x) = \frac{1}{v_{\rm cl}(x)}\int_{0}^{10D} v_x(x, y)\, dy,$$ and the spread $S(x) = \frac{dL}{dx}$.
This one dimensional profiles along the $x$-axis shows that, as expected, the centerline velocity $v_{\rm cl}(x)$ decays proportional to $1/\sqrt{x}$, and that the jet width scales linearly with $x$. In particular, for this simulation the average growth rate (spread) is about $0.1 D$.

\begin{figure}
\includegraphics[width=0.45\textwidth]{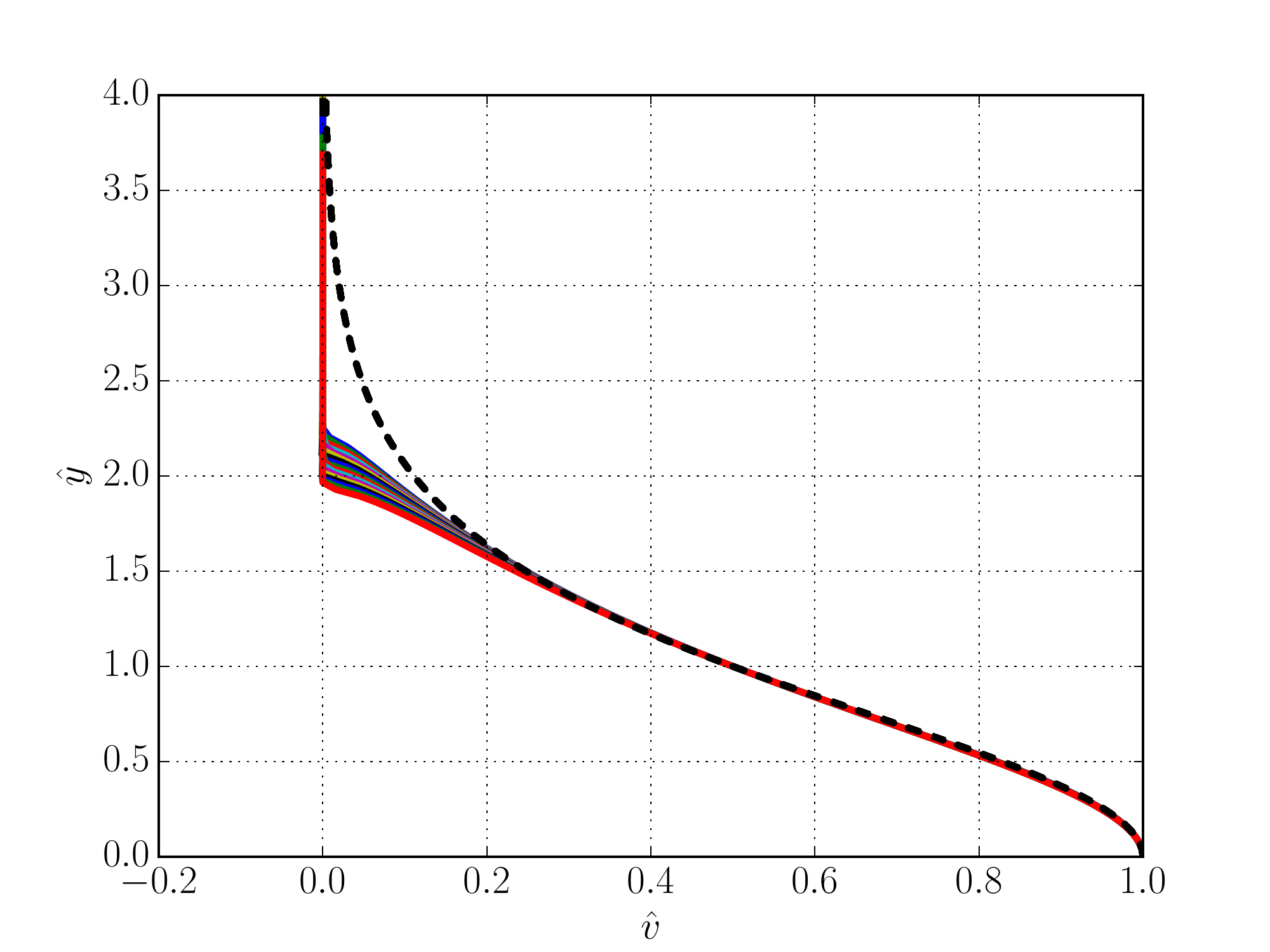}
\includegraphics[width=0.45\textwidth]{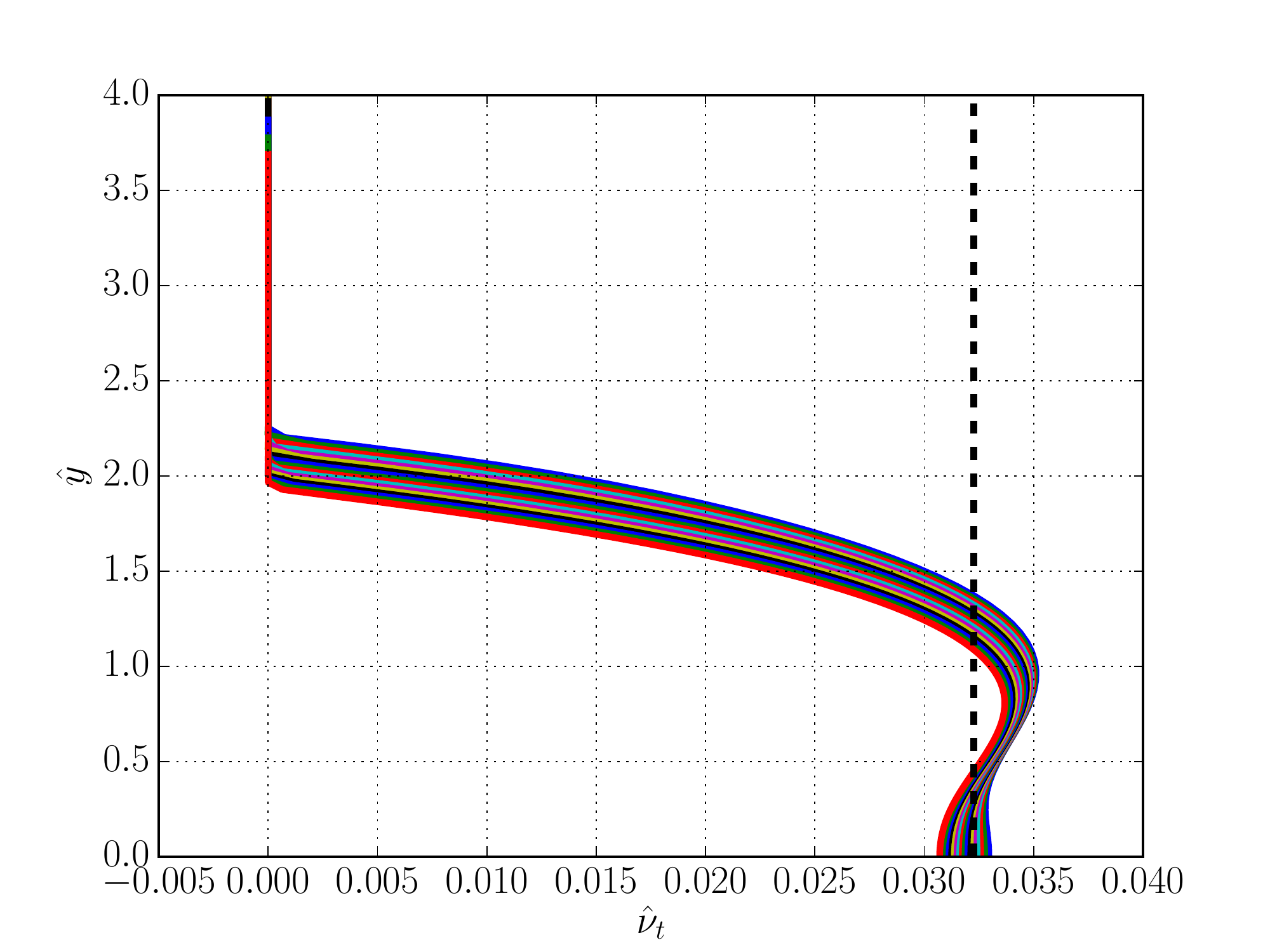}
\caption{Adimensional velocity (left) and turbulent viscosity (profiles) at a distance ranging from $15D$ to $25D$. Dotted lines represent analytical values.}
\label{fig:adim_prof}
\end{figure}

In Figure \ref{fig:adim_prof}, we show adimensional profile for the horizontal velocity and turbulent viscosity at different vertical section at selected distances from the inlet. To this aim, we define the adimensional coordinates $\hat{x} = \frac{x}{D}$ and $\hat{y} = \frac{y}{y_{1/2}}$, the adimensional velocity $\hat{v} = \frac{v_x}{v_{\rm cl}}$, and the adimensional turbulent viscosity $\hat{ \nu_t } = \frac{\nu_t}{ v_{\rm cl}\,y_{1/2}}$.
Figure \ref{fig:adim_prof} shows the profile of $\hat{v}(\hat{ y })$ and $\hat{ \nu_t }(\hat{ y })$. Note how all the adimensional profiles taken at distances ranging from $15D$ to $25D$ from the jet nozzle tends to collapse on the same line. 
The dotted black line in Figure \ref{fig:adim_prof} (left) represents $\hat{v} = \operatorname{sech}^2(\alpha\hat{y})$ with $\alpha = \ln(1 + \sqrt{2})$ which is the analytical velocity profile derived under the assumption of constant turbulent viscosity, see \cite{pope2000turbulent}.
The dotted black vertical line in Figure \ref{fig:adim_prof} (right) represent the reference value $\hat{\nu_t} = 1/31$, see \cite{pope2000turbulent}.
Note how close to the center of the jet the adimensional profiles computed with our finite element solver agree with the reference profiles for this type of flow.

\section{Conclusions}
In this report we discuss a finite element formulation of a $k$-$\epsilon$
closure model for RANS simulations of turbulent free plane jet in two dimension.
The main challenge was to construct a UQ-ready flow solver, i.e., one that is smoothly differentiable and provides not only forward solver capabilities but also adjoint and higher derivatives capabilities, which are needed for scalable Hessian-based UQ methods.
The finite element discretization was implemented in FEniCS and a damped Newton
method was used to solve the nonlinear system of equations arising from
discretization of the steady-state RANS equations.
Velocity and turbulent viscosity profiles computed with our solver
were successfully validated against those from the literature. Also
the centerline velocity, jet width, and growth rate showed the
expected scaling as a function of distance from the inlet boundary.

\bibliographystyle{ieeetr}
\bibliography{plane_jet.bib}

\end{document}